\documentclass[final]{article}

\usepackage{amssymb,amsmath,amsthm,color,graphicx,subfigure}
\usepackage{booktabs}
\usepackage{siunitx}
\usepackage{epstopdf}
\renewcommand{\L}{\mathbf{L}}

\newcommand{\I}{\mathbf{I}}

\newif\iflong
\longtrue % use \longfalse for the publishable paper, \longtrue for the technical report

\newtheorem{theorem}{Theorem}[section]
\newtheorem{remark}{Remark}[section]
\newtheorem{example}{Example}[section]

\usepackage{url}

\title{Application of approximate matrix factorization to high order linearly implicit Runge-Kutta methods}
\author{Hong Zhang\thanks
{Computational Science Laboratory, Department of Computer Science, Virginia Polytechnic Institute and State University, Blacksburg, VA 24061 (zhang@vt.edu)}
\and 
Adrian Sandu\thanks
{Computational Science Laboratory, Department of Computer Science, Virginia Polytechnic Institute and State University, Blacksburg, VA 24061 (sandu@cs.vt.edu)}
\and
Paul Tranquilli\thanks
{Computational Science Laboratory, Department of Computer Science, Virginia Polytechnic Institute and State University, Blacksburg, VA 24061 (ptranq@vt.edu)}
}
\date{}

\begin{document}
\thispagestyle{empty}
\setcounter{page}{0}

\begin{Huge}
\begin{center}
Computer Science Technical Report CSTR-{14-10} \\
\today
\end{center}
\end{Huge}
\vfil
\begin{large}
\begin{center}
Hong Zhang and Adrian Sandu and Paul Tranquilli
\end{center}
\end{large}

\vfil
\begin{huge}
\begin{it}
\begin{center}
Application of approximate matrix factorization to high order linearly implicit Runge-Kutta methods
\end{center}
\end{it}
\end{huge}
\vfil

\begin{large}
\begin{center}
Computational Science Laboratory \\
Computer Science Department \\
Virginia Polytechnic Institute and State University \\
Blacksburg, VA 24060 \\
Phone: (540)-231-2193 \\
Fax: (540)-231-6075 \\ 
Email: \url{sandu@cs.vt.edu} \\
Web: \url{http://csl.cs.vt.edu}
\end{center}
\end{large}

\vspace*{1cm}

\begin{tabular}{ccc}
\includegraphics[width=2.5in]{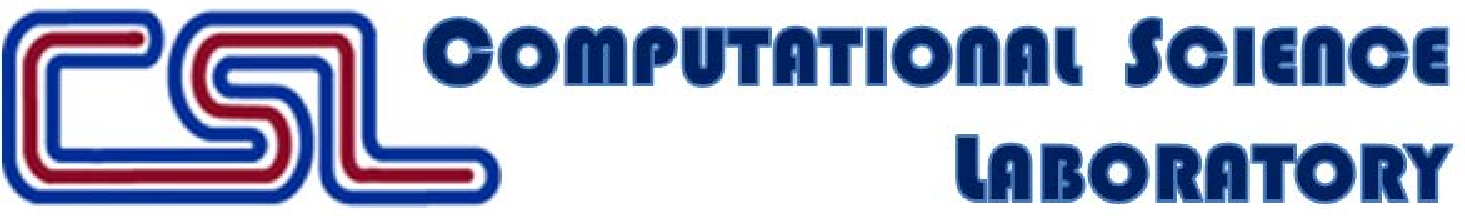}
&\hspace{2.5in}&
\includegraphics[width=2.5in]{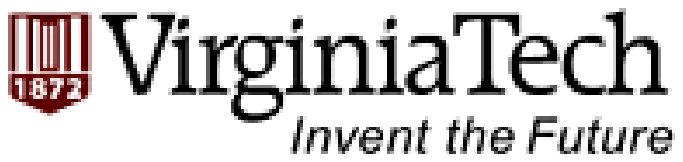} \\
{\bf\em Innovative Computational Solutions} &&\\
\end{tabular}

\newpage

\maketitle

\begin{abstract}
Linearly implicit Runge-Kutta methods with approximate matrix factorization can solve efficiently large systems of differential equations that have a stiff linear part, e.g. reaction-diffusion systems. However, the use of approximate factorization usually leads to loss of accuracy, which makes it attractive only for low order time integration schemes. This paper discusses the application of approximate matrix factorization with high order  methods; an inexpensive correction procedure applied to each stage allows to retain the high order of the underlying linearly implicit Runge-Kutta scheme.  The accuracy and stability of the methods are studied. 
Numerical experiments on reaction-diffusion type problems of different sizes and with different degrees of stiffness illustrate the efficiency of the proposed approach. 
\end{abstract}

%%%%%%%%%%%%%%%%%%%%%%%%%%%%%%%%%%
\section{Introduction}
%%%%%%%%%%%%%%%%%%%%%%%%%%%%%%%%%%
A frequently used approach to solve partial differential equations (PDEs) is the method of lines, 
where the spatial derivative terms are discretized first using techniques such as finite differences, finite volumes or finite elements, and then integrating the resulting system of ordinary differential equations (ODEs) in time. 
Discretization of PDEs with linear terms in space leads to a semi-linear ODE system of the form
\begin{equation}
\label{eqn:semi-linear}
y' = F(y,t) = \L\, y +  f(y,t)\,, \quad y \in \mathbb{R}^N
\end{equation}
where $\L$ is a spatial linear operator and $f(y,t)$ is nonlinear. We consider the case where the linear term has a fast characteristic time scale and the nonlinear term has a slow characteristic time scale.  Due to the Courant-Friedrichs-Levy (CFL) condition, the step size of an explicit time integrator is restricted by the fastest time scale. 
%The problem is stiff when the step size is forced to values that are much smaller than necessary to obtain the desired accuracy. 
To alleviate the time step constraint imposed by the stiff linear term with a reasonably low computational cost, 
linearly implicit methods (particular cases of implicit-explicit methods) treat the stiff linear term implicitly while the nonlinear term explicitly.  
%
% The semi-linear partitioning can be the result of an artificial splitting of the right hand side: \textcolor{red}{Do we need this splitting in this paper?}
% %
% \begin{equation}
% \label{eqn:splitting}
% y' = F(y) = \L\, y + \left( F(y) - \L\,y \right) \, = g(y) + f(y)\,, \quad y \in \mathbb{R}^N
% \end{equation}

Various families of linearly implicit integration methods have been proposed  and successfully applied to solve PDEs with linear dispersion and dissipation \cite{Grooms_2011,Calvo_2001,Akrivis_2003_IMEX,Akrivis_2004}.
Fully implicit schemes solve at each step a linear system where the matrix involves the Jacobian of the right hand side function. 
Efficient schemes specialized in the solution of \eqref{eqn:semi-linear} such as linearly implicit methods use only the part of the Jacobian associated with the stiff linear term, 
resulting in linear systems of the form
\begin{equation}
\label{eqn:linearsys}
\left( \I - h\, \gamma \, \L\right)\,x = \ell, \quad x,\ell \in \mathbb{R}^N,
\end{equation}
where $h$ is the step size, $\gamma$ is a parameter of the integration scheme, and the right hand side $\ell$ is determined by the method. 

The matrix $\L$  is usually sparse but has a large bandwidth, especially when high order spatial discretization schemes are applied. 
Consequently the LU factorization of the matrix in \eqref{eqn:linearsys} can be very costly in large scale problems. 
One approach to increase efficiency is to apply iterative solvers to \eqref{eqn:linearsys}, as in \cite{Weiner_1997_rowmap}, 
%where Rosenbrock-type methods like ROW- and W-mehtods are applied for stiff problems, 
but there are associated challenges related to preconditioning and convergence. 

An alternative approach to increase the computational efficiency is approximate matrix factorization (AMF).
Assuming that the matrix is a sum of simpler matrices
\begin{equation}
\label{eqn:LM}
  \L  = \sum_{r=1}^R \L_r
\end{equation}
AMF replaces the system matrix \eqref{eqn:linearsys} with a product of simpler, and easier to factorize, matrices
\begin{equation}
\label{eqn:amf}
 \I - h\, \gamma \, \L  \approx  \I - h\, \gamma \, \widetilde{\L} = \prod_{r=1}^R \left(\I - h\, \gamma \, \L_r\right).
\end{equation}
The approximation formula \eqref{eqn:amf} defines implicitly the matrix $\widetilde{\L}$ as 
\begin{equation}
\label{eqn:L}
 \widetilde{\L}  = \L + \sum_{k=2}^R \, \left(-h \gamma\right)^{k-1} \sum_{1\leq i_1 < i_2 < \dots < i_k \leq R} \L_{i_1} \, \L_{i_2} \dots \L_{i_k} .
 % h\, \gamma \sum_{i=1}^{R-1} \sum_{j=i+1}^{R} \L_i \, \L_j 
 % +  h^2\, \gamma^2 \sum_{i=1}^{R-2} \sum_{j=i+1}^{R-1} \sum_{k=j+1}^{R} \L_i \, \L_j \, \L_k  + \dots  + (-1)^R h^R \gamma^R \prod_{i=1}^R \L_i.
 % \mathcal{O}(h^2). 
\end{equation}
For example, consider $\L$ to be the discrete two-dimensional Laplace operator.  It can be written in the form \eqref{eqn:LM} where $\L_1$ and $\L_2$ correspond to the derivatives along the x direction and y direction, respectively. 
The AMF approximation corresponds to the alternating directions factorization \cite{Douglas_1962,Peaceman_1955}
\[
\I - h\, \gamma \, \widetilde{\L} := \left(\I-h\, \gamma \L_1\right) \, \left(\I -h\, \gamma \L_2\right)
\]
where 
\begin{equation}
\label{eqn:L2}
\widetilde{\L}  = \L -h\, \gamma\, \L_1 \, \L_2. 
\end{equation}

The idea of using AMF to speed up calculations in implicit time integration has appeared multiple times in the literature.
Sandu \cite{Sandu_1997_phd} discussed a family of methods named ELADI, which are Rosenbrock-W schemes that make use of AMF to speed up calculations. 
The order of the resulting discretization remains unchanged since Rosenbrock-W schemes can accommodate any Jacobian approximation.
Houwen et al. provided a survey for AMF methods applied in the context of several different linear integration schemes, and provide stability results for such schemes \cite{Houwen_2001}.  The discussion is limited to methods of low order, two and three, presumably due to the inaccurate nature of the AMF approximation.  
Gonzalez \cite{Gonzalez_2014} proposed a way to apply AMF-refinements to second-order and third-order Rosenbrock-type methods for solving advection-diffusion-reaction PDEs. 
But their methods are limited to the low or medium precision level, and no generalization to higher order is supplied.
In \cite{Beck_2014} Beck et al. compared the efficiency of AMF versus Krylov based approaches to the solution of linear systems in the context of Newton iterations arising in Radau \cite{Hairer_book_II} and Peer \cite{Beck_2012} integration methods.  These methods avoid the issue of order degradation, through the use of integration schemes in which the Jacobian of the spatial discretization does not appear explicitly. 
They conclude that AMF methods are extremely efficient, particularly when low accuracy solutions are sought.  
Berzins et al. \cite{Berzins_1996, Ahmad_1997} presented a method for solving the linear system in a Newton iteration arising from several classes of time integration methods including theta methods, backward differential formulas, and implicit-explicit (IMEX) multistep methods.
They performed and analysis of the error arising from operator splitting and provided a method to control timesteps such as to guarantee Newton convergence when using AMF. 
Since AMF is only used to speed up the solution of the nonlinear equations, the error does not affect the order of accuracy at the time stepping level.   
%However, existing AMF approaches limit the application only to low order schemes due to the possible loss of order resulting from the AMF approximation \cite{Calvo_2005_IMEX_AMF} . 
%Our goal is to apply the AMF technique to high order methods, at the same time maintaining the high order of accuracy and high efficiency.     
%
Calvo and Gerisch \cite{Calvo_2005_IMEX_AMF} applied AMF to a form of linearly implicit Runge-Kutta (LIRK) methods that avoids the computation of matrix-vector products. 
First order of convergence is obtained by using a third-order LIRK method, and improved to second order by adding a correction to the solution at each time step. 

None of the existing methods that take advantage of AMF can provide highly accurate results and there is still room for improvement in efficiency. 
The goal of this work is to achieve high accuracy while maintaining a low computational cost. 
%Runge-kutta type methods are known to have good stability properties. 
%The stability regions of these methods usually cover a large part of the imaginary axis.
%This is essential for solving problems that have complex eigenvalues which are not uncommon in advection-diffusion-reaction PDEs. 
The focus is on using AMF with linearly implicit Runge-Kutta methods of high order. 
The main contribution of this paper is to account for the inherent inaccuracy of AMF through low cost refinements of the stage values and to recover the accuracy of the underlying time discretization. 
%The theoretical findings are validated  with numerical experiments. 

The remainder of this paper is organized as follows. 
Section \ref{sec:LIRK} introduces LIRK methods and the existing approaches to apply AMF. 
Section \ref{sec:stage_refinement} presents a new strategy to incorporate AMF by using an inexpensive stage refinement procedure. 
An error analysis explains how this strategy solves the accuracy degradation issue that affects existing approaches.  Stability issues are also investigated. 
Section \ref{sec:AMF_exp} reports numerical rests for a variety of test problems of different dimensions and different degrees of stiffness, and illustrates the convergence behavior and efficiency of the approach.  Conclusions are drawn in Section \ref{sec:AMF_con}.

%%%%%%%%%%%%%%%%%%%%%%%%%%%%%%%%%%%%
\section{Linearly implicit Runge-Kutta methods} \label{sec:LIRK}
%%%%%%%%%%%%%%%%%%%%%%%%%%%%%%%%%%%%

A general linearly implicit Runge-Kutta (LIRK) scheme proposed by Calvo, de Frutos, and Novo \cite{Calvo_2001} is obtained 
by applying the IMEX Runge-Kutta methods \cite{Ascher_1997,Pareschi_2000}
\begin{subequations}
\label{eqn:imex_rk}
\begin{eqnarray}
Y_i  &=& y_n +   h\, \sum_{j=1}^{i-1}\, a_{i,j}\,   f(Y_j)  +   h\, \sum_{j=1}^i\, \widehat{a}_{i,j}\,   g(Y_j), \\
y_{n+1} &=&  y_n +   h\, \sum_{j=1}^s\, b_{j}\,   f(Y_j)  +   h\, \sum_{j=1}^s\, \widehat{b}_{j}\,   g(Y_j), 
\end{eqnarray}
\end{subequations}
to solve \eqref{eqn:semi-linear} where the stiff component is linear, $g(y) := \L y$:
\begin{subequations}
\label{eqn:lirk}
\begin{eqnarray}
\label{eqn:lirk-stage}
\left( \I - h\,\widehat{a}_{i,i}\,\L\right)\, Y_i  &=& \ell_i := y_n +   h\, \sum_{j=1}^{i-1}\, a_{i,j}\,   f(Y_j)  +  h\, \sum_{j=1}^{i-1}\,  \widehat{a}_{i,j} \,  \L  \, Y_j, \\
\label{eqn:lirk-step}
y_{n+1} &=&  y_n +   h\, \sum_{j=1}^s\, b_{j}\,   f(Y_j)  +   h\, \, \sum_{j=1}^s\,  \widehat{b}_{j} \L \,  Y_j. 
\end{eqnarray}
\end{subequations}
Order conditions for LIRK methods are derived in \cite{Calvo_2001}, and simplifications of the IMEX Runge-Kutta conditions are possible due to the special form of the nonlinear term.

The coefficients $\widehat{a}_{i,i}$ in practical methods are chosen to be equal to the same value $\gamma$ for computational efficiency, as this allows to reuse the same LU factorization in all stages.
%
%Moreover the stage equations \eqref{eqn:lirk-stage} can be written as
%%
%\begin{eqnarray*}
%\left( \I - h\,\widehat{a}_{i,i}\,\L\right)\, Y_i  &=& y_n +   h\, \sum_{j=1}^{i-1}\, a_{i,j}\,   F(Y_j)  +  h\, \L  \, \sum_{j=1}^{i-1}\,  (\widehat{a}_{i,j}-a_{i,j}) \,  Y_j 
%\end{eqnarray*}
%
A linear transformation of variables allows for a reformulation of the stage equations \eqref{eqn:lirk-stage} in a form that avoids explicit multiplications by the matrix $\L$
\begin{subequations}
\label{eqn:lirk-no-multiplication}
\begin{eqnarray}
\label{eqn:lirk-no-multiplication-stage-u}
\left( \I - h\,\gamma\,\L\right)\,U_i  &=& y_n +  \sum_{j=1}^{i-1}\,  \frac{ \widehat{a}_{i,j}-a_{i,j} }{\gamma}\,   \, Y_j + h\, \sum_{j=1}^{i-1}\, a_{i,j}\,   F(Y_j), \\
\label{eqn:lirk-no-multiplication-stage}
Y_i &=& U_i -  \sum_{j=1}^{i-1}\,  \frac{ \widehat{a}_{i,j}-a_{i,j} }{\gamma}\,   \, Y_j. 
\end{eqnarray}
Moreover, it is convenient to choose pairs of methods with the same weights $b_{j} = \widehat{b}_{j}$ in which case
the next step solution \eqref{eqn:lirk-step} is 
\begin{equation}
\label{eqn:lirk-no-multiplication-step}
y_{n+1} =  y_n +   h\, \sum_{j=1}^s\, b_{j}\,  F(Y_j). 
\end{equation}
\end{subequations}

%%%%%%%%%%%%%%%%%%%%%%%%%%%%
%\section{Approximate matrix factorization} \label{sec:amf}
%%%%%%%%%%%%%%%%%%%%%%%%%%%%

Calvo and Gerisch \cite{Calvo_2005_IMEX_AMF} studied the use of LIRK methods with AMF. The approximation is obtained by replacing $\I - h \gamma \L$ with $\I - h \gamma \widetilde{\L}$ in \eqref{eqn:lirk-no-multiplication-stage-u}, or equivalently, by using the matrix $\widetilde{\L}$ instead of $\L$ in \eqref{eqn:lirk}
\begin{eqnarray*}
\left( \I - h\,\gamma \, \widetilde{\L} \right)\, \widetilde{Y}_i  &=& y_n +   h\, \sum_{j=1}^{i-1}\, a_{i,j}\,   f(\widetilde{Y}_j)  +  h\, \sum_{j=1}^{i-1}\,  \widehat{a}_{i,j} \, \widetilde{\L}  \, \widetilde{Y}_j, \\
\widetilde{y}_{n+1} &=&  y_n +   h\, \sum_{j=1}^s\, b_{j}\,   f(\widetilde{Y}_j)  +   h\, \sum_{j=1}^s\, \widehat{b}_{j} \,   \widetilde{\L} \,\widetilde{Y}_j.
\end{eqnarray*}
It can also be regarded as a direct application of the LIRK method \eqref{eqn:lirk} to solve the perturbed ODE system 
\begin{equation}
\label{eqn:splitting-new}
%y' = \widetilde{\L} \, y + f(y) = F(y) - h\, \gamma\, \left( \sum_{i=1}^{M-1} \sum_{j=i+1}^{M} \L_i \, \L_j  \right) \, y + \mathcal{O}(h^2)
\widetilde{y}' = \widetilde{\L} \, \widetilde{y} + f(\widetilde{y}) = F(\widetilde{y}) - h\, \gamma\, \left( \L - \widetilde{\L} \right) \, \widetilde{y}
\end{equation}
instead of the original system \eqref{eqn:semi-linear}.
The first-order behavior of this approach has been explained in \cite{Calvo_2005_IMEX_AMF} by the fact that the perturbation added in \eqref{eqn:splitting-new} 
changes the solution over one time step by $\mathcal{O}(h^2)$.

To recover second order Calvo and Gerisch \cite{Calvo_2005_IMEX_AMF} apply corrections to the numerical solution obtained by LIRK with AMF. One such correction has the form
\begin{equation}
\label{eqn:calvo-correction}
y_{n+1} = \widetilde{y}_{n+1} 
%+ h^2\, \gamma\,  \left( \I - h\,\gamma \, \widetilde{\L} \right)^{-1}\, \left( \sum_{i=1}^{M-1} \sum_{j=i+1}^{M} \L_i \, \L_j \right)\, y_n\,.
+ h^2\, \gamma\,  \left( \I - h\,\gamma \, \widetilde{\L} \right)^{-1}\, \left( \L - \widetilde{\L} \right)\, y_n\,.
\end{equation}
The matrix inverse in the correction term uses the same LU factorization as the solution. The presence of the matrix inverse in the correction term ensures the linear stability of the new solution \eqref{eqn:calvo-correction}.
It is also noted in \cite{Calvo_2005_IMEX_AMF} that ``regaining order three using corrections similar to \eqref{eqn:calvo-correction} is not feasible with a computational cost comparable with the cost of recovering order two''.

%Since the linear term is treated implicitly, the perturbation on $\L$ will not affect the integration scheme in terms of stability. 
%But it brings an extra error of the magnitude $O(h\gamma \|\L_1\,\L_2\|)$ to the right hand side, which is the first derivative of the ODE solution. 
% The analysis in \cite{Calvo_2005_IMEX_AMF} reveals that the local error for one-step integration without corrections is $\mathcal{O}\left(h^2\gamma \|\L_1\,\L_2 \|y_n \right)$. Need to revisit this - a local error of order $h^2$ means order one, NOT order two. Please check!

%%%%%%%%%%%%%%%%%%%%%%%%%%%
\section{LIRK-AMF methods with stage refinement} \label{sec:stage_refinement}
%%%%%%%%%%%%%%%%%%%%%%%%%%%

We consider a different way to incorporate AMF into LIRK methods, which forms the basis of all approaches presented in this paper. 
In order to keep the right-hand side of the original ODE system \eqref{eqn:semi-linear}, 
we only approximate $ \I - h\,\gamma \L$ with $\I - h\,\gamma \widetilde{\L}$ when computing the Runge-Kutta stages: 
\begin{subequations}
\label{eqn:lirk-hong}
\begin{eqnarray}
\label{eqn:lirk-hong-stage}
\left( \I - h\, \gamma \,\widetilde{\L}\right)\, \widetilde{Y}_i  &=& \widetilde{\ell_i} := y_n +   h\, \sum_{j=1}^{i-1}\, a_{i,j}\,   f(\widetilde{Y}_j)  +  h\, \sum_{j=1}^{i-1}\,  \widehat{a}_{i,j} \,  \L  \, \widetilde{Y}_j, \\
\label{eqn:lirk-hong-step}
\widetilde{y}_{n+1} &=&  y_n +   h\, \sum_{j=1}^s\, b_{j}\,   f(\widetilde{Y}_j)  +   h\, \sum_{j=1}^s\,  \widehat{b}_{j} \L \,  \widetilde{Y}_j. 
\end{eqnarray}
\end{subequations}
Thus the new method uses an inexact Jacobian for the implicit part in the LIRK scheme. 
The change of the left hand side in \eqref{eqn:lirk-hong-stage} will affect the solution, and a correction is needed to restore accuracy.

Consider the solution of the original stage equation \eqref{eqn:lirk-stage} 
\[
\left( \I - h \, \gamma \, \L \right)\, Y_i - \ell_i = 0
\]
by simplified Newton iterations of the form:
\begin{equation}
\label{eqn:simplified-Newton}
 Y_i^{(k)} =  Y_i^{(k-1)} - \left( \I - h \, \gamma \, \widetilde{\L} \right)^{-1} \cdot \left( \left( \I - h \, \gamma \, \L \right)\, Y_i^{(k-1)} - \ell_i \right).
\end{equation}
For example, the direct solution is:
\begin{eqnarray*}
Y_i^{(0)} &=& \left( \I - h \, \gamma \, \widetilde{\L} \right)^{-1}\, \ell_i,
\end{eqnarray*}
and the solution after one refinement iteration is:
\begin{eqnarray}
\label{eqn:one-ieration}
 Y_i^{(1)} &=&  Y_i^{(0)} - \left( \I - h \, \gamma \, \widetilde{\L} \right)^{-1} \cdot \left( \left( \I - h \, \gamma \, \L \right)\, Y_i^{(0)} - \ell_i \right)\\
 \nonumber
&=& Y_i^{(0)} -  \left( \I - h \, \gamma \, \widetilde{\L} \right)^{-1} \cdot \left( h \, \gamma \,\widetilde{\L} - h \, \gamma \,\L \right)\, Y_i^{(0)}.
% \\
%&=& \left( \I - h \, \gamma \, \widetilde{\L} \right)^{-1} \cdot \left( \I - 2\, h \, \gamma \, \widetilde{\L}  + h \, \gamma \,\L \right)\, Y_i^{(0)}.
\end{eqnarray}

Next we analyze the linear system solution errors and  investigate how this errors propagate to affect the solution at the next step. 
\subsection{Error analysis} 
%%%%%%%%%%%%%%%%%%%%%%

Consider the exact stage solution
\[
 Y_i = \left( \I - h \, \gamma \, \L \right)^{-1}\,\ell_i. 
\]
The linear system solution error after $k$ iterations is defined as 
\[
 \varepsilon_i^{(k)}  =  Y_i^{(k)}  -  Y_i.
 \]
From \eqref{eqn:simplified-Newton} we obtain
\begin{subequations}
\label{eqn:error-equation}
\begin{eqnarray}
\varepsilon_i^{(k)} &=&  \varepsilon_i^{(k-1)} - \left( \I - h \, \gamma \, \widetilde{\L} \right)^{-1} \cdot  \left( \I - h \, \gamma \, \L \right)\, \varepsilon_i^{(k-1)} \\
 &=&   \left( \I - h \, \gamma \, \widetilde{\L} \right)^{-1} \cdot  \left(  \I - h \, \gamma \, \widetilde{\L}  - (\I - h \, \gamma \, \L) \right)\, \varepsilon_i^{(k-1)} \\
 \label{eqn:error-decrease}
 &=& -  \left( \I - h \, \gamma \, \widetilde{\L} \right)^{-1}  \cdot  \left( h\, \gamma\,\widetilde{\L}  - h\, \gamma\,\L \right)\, \varepsilon_i^{(k-1)}. 
% \\
% &=& -h\, \gamma\,  \left( \I - h \, \gamma \, \widetilde{\L} \right)^{-1}  \cdot  \left(  -h\, \gamma\,  \sum_{i=1}^{M-1} \sum_{j=i+1}^{M} \L_i \, \L_j + \dots \right) \, \varepsilon_i^{(k-1)}. 
\end{eqnarray}
\end{subequations}

\paragraph{Nonstiff or moderately stiff case.} In the nonstiff or moderately stiff case we have $\Vert h \L_i \Vert = \mathcal{O}(h)$, therefore $\Vert \widetilde{\L} -\L \Vert = \mathcal{O}(h)$ and 
\[
\left\Vert  \left( \I - h \, \gamma \, \widetilde{\L} \right)^{-1}  \cdot  \left( h \, \gamma \,\widetilde{\L}  - h \, \gamma \,\L \right) \right\Vert = \mathcal{O}\left(h^2\right).
\]
Consequently from  \eqref{eqn:error-decrease} the error decrease is
\begin{eqnarray*}
\left\Vert \varepsilon_i^{(k+1)} \right\Vert &=&   \mathcal{O}(h^2) \, \left\Vert \varepsilon_i^{(k)} \right\Vert \quad \Rightarrow \quad  \left\Vert \varepsilon_i^{(k)} \right\Vert = \mathcal{O}\left(h^{2k+2}\right).
\end{eqnarray*}

\paragraph{Highly stiff case.} For the highly stiff case $\Vert h \L_i \Vert \gg 1$. We make the assumption
that, for any $\overline{h}$ there exists $0 < \rho(\overline{h}) < 1$ such that the following matrix norm is uniformly bounded for any step size smaller than $\overline{h}$:
\[
\left\Vert  \left( \I - h \, \gamma \, \widetilde{\L} \right)^{-1}  \cdot  \left( h \, \gamma \,\widetilde{\L}  - h \, \gamma \,\L \right) \right\Vert \le \rho(\overline{h}) < 1\,, \quad \forall h:~ 0 \le h \le \overline{h}.
\]
In the highly stiff case  the error decrease equation \eqref{eqn:error-decrease} leads to
\begin{eqnarray*}
 \left\Vert \varepsilon_i^{(k)} \right\Vert &=&   \rho \, \left\Vert \varepsilon_i^{(k-1)} \right\Vert \quad \Rightarrow \quad \left\Vert \varepsilon_i^{(k)} \right\Vert = \rho^{k}\, \varepsilon_i^{(0)} .
\end{eqnarray*}
We expect that the convergence rate will decrease with increasing step sizes, i.e., $\rho(\overline{h}) \to 1$ when $\overline{h} \to \infty$.

\begin{example}[Dimensional splitting of the discrete diffusion operator on a Cartesian grid] Consider the two-dimensional diffusion operator with periodic boundary conditions on a domain of size $L_X \times L_Y$. 
It is discretized on an $M \times N$ grid of size $\Delta x$, $\Delta y$.
We perform a dimensional splitting.
The error equation  \eqref{eqn:error-decrease}
%
%\[
%\left( \I - h \, \gamma \, \widetilde{\L} \right)\, \varepsilon_i^{(k+1)} = -  \left( h\, \gamma\,\widetilde{\L}  - h\, \gamma\,\L \right)\, \varepsilon_i^{(k)} 
%\]
%
can be written as
\begin{equation}
\label{eqn:error-eqn-diffusion}
\left( \I - h \, \gamma \, \L_1 \right)\, \left( \I - h \, \gamma \, \L_2 \right)\,\varepsilon_i^{(k)} = \left( h\, \gamma\,\L_1  \right)\, \left( h\, \gamma\,\L_2  \right)\,\varepsilon_i^{(k-1)}. 
\end{equation}
Consider the discrete frequencies
\[
-\frac{M}{2} \le m \le \frac{M}{2}-1, \quad -\frac{N}{2} \le n \le \frac{N}{2}-1, \quad 
\tilde{m} = \frac{2\pi m}{L_X},% = \frac{2\pi m}{M\, \Delta x} 
 \quad 
\tilde{n} = \frac{2\pi n}{L_Y}. % = \frac{2\pi n}{N\, \Delta y}
\]
A discrete Fourier transform applied to \eqref{eqn:error-eqn-diffusion} gives the following error equation for each spatial mode $(m,n)$ of the error: 
\[
(1 + h \, \gamma \, \tilde{m}^2)(1 + h \, \gamma \, \tilde{n}^2)\, \hat{\varepsilon}^{(k)}_{m,n} =  
(- h \, \gamma \, \tilde{m}^2)(- h \, \gamma \, \tilde{n}^2)\, \hat{\varepsilon}^{(k-1)}_{m,n}
\]
Let
\[
z_1 = \frac{h}{\Delta x^2}\, \left( \frac{2\pi m}{M} \right)^2, \quad
z_2 = \frac{h}{\Delta y^2}\, \left( \frac{2\pi n}{N} \right)^2.
\]
The evolution of the mode $(m,n)$ of the error is
\[
\hat{\varepsilon}^{(k+1)}_{m,n} =  \frac{(\gamma z_1)(\gamma z_2)}{(1+\gamma z_1)(1+\gamma z_2)} \hat{\varepsilon}^{(k)}_{m,n}
\]
The error amplification factor for the $(m,n)$ mode is
\[
R_{m,n} = \frac{(\gamma z_1)(\gamma z_2)}{(1+\gamma z_1)(1+\gamma z_2)}, \quad
| R_{m,n} | < 1, \quad
| R_{m,n} | \stackrel{z_1,z_2 \to \infty}{\longrightarrow} 1. 
\]
Therefore we expect that more iterations will be required for stiff problems. The AMF with correction will work well for mildly stiff problems. It will work well for stiff problems only when the solution is smooth, and the high order modes are approaching zero.
Similar conclusions are drawn for the three-way splitting of a three dimensional diffusion problem
\[
R_{m,n.k} = \frac{(\gamma z_1)(\gamma z_2)+(\gamma z_1)(\gamma z_3)+(\gamma z_2)(\gamma z_3)+(\gamma z_1)(\gamma z_2)(\gamma z_3)}{(1+\gamma z_1)(1+\gamma z_2)(1+\gamma z_3)}.
\]
\end{example}

\begin{remark}
The accuracy analysis in Calvo and Gerisch's paper \cite{Calvo_2005_IMEX_AMF} considers the non-stiff case.
For very stiff systems the correction term \eqref{eqn:calvo-correction} is
\[ 
h\, \left( \I - h\,\gamma \, \widetilde{\L} \right)^{-1}\, \left( h\, \gamma\,\L - h\, \gamma\,\widetilde{\L} \right)\, y_n
= \mathcal{O}(h)
\]
and the remaining error term is $\mathcal{O}(h^2)$, therefore the corrected solution is first order.
\end{remark}

%%%%%%%%%%%%%%%%%%%%%%%%%%%%%%%%%%%%%%%%%
\subsection{Propagation of linear system errors} 
%%%%%%%%%%%%%%%%%%%%%%%%%%%%%%%%%%%%%%%%%
The computation of stage values via \eqref{eqn:lirk-hong} and \eqref{eqn:simplified-Newton}  propagates the linear system errors from one stage to another.  To account for the total error consider the methods \eqref{eqn:lirk} and \eqref{eqn:lirk-hong} and let
\begin{equation}
\label{eqn:def_delta_Y}
\delta {Y}_i  =  \widetilde{Y_i}  -  {Y}_i.
\end{equation}
We assume that these errors are small.
The exact stage equations \eqref{eqn:lirk-stage} read
\begin{eqnarray*}
\left( \I - h\,\gamma \,{\L}\right)\, {Y}_i  &=& y_n +   h\, \sum_{j=1}^{i-1}\, a_{i,j}\,   f({Y}_j)  +  h\, \sum_{j=1}^{i-1}\,  \widehat{a}_{i,j} \,  \L  \,{Y}_j.
%&=& y_n +   h\, \sum_{j=1}^{i-1}\, a_{i,j}\,   f(\widetilde{Y}_j-\delta {Y}_j)  +  h\, \sum_{j=1}^{i-1}\,  \widehat{a}_{i,j} \,  \L  \,(\widetilde{Y}_j-\delta {Y}_j) \\ 
%&=& y_n +   h\, \sum_{j=1}^{i-1}\, a_{i,j}\,   f(\widetilde{Y}_j)  +  h\, \sum_{j=1}^{i-1}\,  \widehat{a}_{i,j} \,  \L  \,(\widetilde{Y}_j) \\ 
%& &   +   h\, \sum_{j=1}^{i-1}\, a_{i,j}\,   f'(\widetilde{Y}_j) \cdot \delta {Y}_j  -  h\, \sum_{j=1}^{i-1}\,  \widehat{a}_{i,j} \,  \L  \,\delta {Y}_j + \mathcal{O}(\Vert \delta {Y} \Vert^2)\\ 
\end{eqnarray*}
The AMF stage equations are solved inexactly and read
\begin{eqnarray*}
\left( \I - h\,\gamma \,\L\right)\, \left( \widetilde{Y}_i - \varepsilon_i \right)  &=& y_n +   h\, \sum_{j=1}^{i-1}\, a_{i,j}\,   f(\widetilde{Y}_j)  +  h\, \sum_{j=1}^{i-1}\,  \widehat{a}_{i,j} \,  \L  \, \widetilde{Y}_j, 
\end{eqnarray*}
where $\varepsilon_i$ is the error due to the simplified Newton approximation of the system solution. Using \eqref{eqn:def_delta_Y} we express this in terms of the solution of the exact stages
\begin{eqnarray*}
\left( \I - h\,\gamma\,\L\right)\, \left( Y_i + \delta Y_i - \varepsilon_i \right)  
\iflong
&=& y_n +   h\, \sum_{j=1}^{i-1}\, a_{i,j}\,   f( Y_j + \delta Y_j )  \\
&& +  h\, \sum_{j=1}^{i-1}\,  \widehat{a}_{i,j} \,  \L  \,( Y_j + \delta Y_j ) \\
\fi
&=& y_n +   h\, \sum_{j=1}^{i-1}\, a_{i,j}\,   f({Y}_j)  +  h\, \sum_{j=1}^{i-1}\,  \widehat{a}_{i,j} \,  \L\, {Y}_j \\ 
& &   +   h\, \sum_{j=1}^{i-1}\, a_{i,j}\,   f'_j \, \delta {Y}_j  +  h\, \sum_{j=1}^{i-1}\,  \widehat{a}_{i,j} \,  \L  \,\delta {Y}_j, 
\end{eqnarray*}
where we use the mean value theorem
\[
 f( Y_j + \delta Y_j ) - f({Y}_j)  = f'_j \cdot \delta {Y}_j\,, \quad f'_j = \int_0^1 f_y(Y_j + s\, \delta Y_j)\, ds.
\]
The nonstiff/moderately stiff assumption about the nonlinear terms $f$ implies that these average Jacobians are of moderate size,
\begin{equation}
\label{eqn:nonstiff-f-jacobian}
\Vert f'_j \Vert = \mathcal{O}(1)\,, \quad \forall\, j\,.
\end{equation}

After subtracting the exact stage equations we are left with the error relation
%
%\begin{eqnarray*}
%\left( \I - h\,\gamma \,\L\right)\, \left(  \delta Y_i - \varepsilon_i \right)  &=&   h\, \sum_{j=1}^{i-1}\, a_{i,j}\,   f'_j \cdot \delta {Y}_j  +  h\, \sum_{j=1}^{i-1}\,  \widehat{a}_{i,j} \,  \L  \,\delta {Y}_j.
%\end{eqnarray*}
%%
%or
%
%\begin{eqnarray}
%\label{eqn:stage-error-equation1}
% \delta Y_i &=&  \varepsilon_i  +  \sum_{j=1}^{i-1}\, a_{i,j}\,  \left( \I - h\,\gamma \,\L\right)^{-1}\, \left(h\, f'_j\right) \cdot \delta {Y}_j  \\
%\nonumber 
%     &&  +  \sum_{j=1}^{i-1}\,  \widehat{a}_{i,j} \,  \left( \I - h\,\gamma \,\L\right)^{-1}\, (h\, \L)  \,\delta {Y}_j.
%\end{eqnarray}
%%
%This can be viewed as
%%
%\begin{eqnarray}
%\label{eqn:stage-error-equation2}
% \delta Y_i &=&  \left( \I - h\,\gamma \,\L\right)\,  \varepsilon_i  +  \sum_{j=1}^{i-1}\, a_{i,j}\,  \left(h\, f'_j\right) \cdot \delta {Y}_j  \\
%\nonumber 
%     &&  +  \sum_{j=1}^{i}\,  \widehat{a}_{i,j} \,  (h\, \L)  \,\delta {Y}_j.
%\end{eqnarray}
%
%In matrix notation with $\mathbf{\widehat{\Gamma}}$=diag$(\mathbf{\widehat{A}})$ we have
%%
%\begin{equation}
%\label{eqn:delta-stages}
%\left( \I - h\,\mathbf{\widehat{A}} \otimes \L - h\,\mathbf{A} \odot F'\ \right)\, \delta Y = \left( \I - h\,\mathbf{\widehat{\Gamma}} \otimes \L\right)\, \varepsilon
%\end{equation}
%%
%or
%
\begin{equation}
\label{eqn:delta-stages-2}
\delta Y = \left( \I - h\,\mathbf{\widehat{A}} \otimes \L - h\,\mathbf{A} \odot F'\ \right)^{-1}\, \left( \I - h\,\mathbf{\widehat{\Gamma}} \otimes \L\right)\, \varepsilon
\end{equation}
where $\mathbf{\widehat{\Gamma}}$=diag$(\mathbf{\widehat{A}})$ and $(\mathbf{A} \odot F')_{i,j} = a_{i,j} f'_j$.
For nonstiff or moderately stiff problems $\Vert h\L \Vert= \mathcal{O}(h)$, $\Vert h F'\Vert = \mathcal{O}(h)$, and for small step sizes we have
\[
\Vert \delta Y \Vert = \left(  1 + \mathcal{O}(h) \right)\, \Vert \varepsilon \Vert.
\]
For highly stiff systems $\Vert h \L \Vert \to \infty$ and we have
\[
\delta Y = \left( \mathbf{\widehat{A}} \otimes (\L/\Vert \L \Vert)  \right)^{-1}\, \left( \mathbf{\widehat{\Gamma}} \otimes (\L/\Vert \L \Vert)\right)\, \varepsilon
\]
therefore
\[
\Vert \delta Y \Vert =  \mathcal{O}(1) \, \Vert \varepsilon \Vert.
\]
In both the nonstiff and the stiff cases the stage error is of the size of the linear system solution error.

%It follows from \eqref{eqn:nonstiff-f-jacobian} that the following matrix norms are uniformly bounded for any step size:
%%
%\[
%\left\Vert\,  \left( \I - h \, \gamma \, \L \right)^{-1}  \cdot  f'_j \, \right\Vert \le C\,, \quad \forall \, h >0.
%\]
%%
%Under these assumptions it is clear that
%%
%\[
%\varepsilon_i \sim \mathcal{O}\left( h^k \right) \quad \Rightarrow \quad   \delta Y_i \sim \mathcal{O}\left( h^k \right)\,.
%\]

From the exact step equation the error in the solution is
%
%\[
%y_{n+1} =  y_n +   h\, \sum_{j=1}^s\, b_{j}\,   f(Y_j)  +   h\, \sum_{j=1}^s\,  \widehat{b}_{j} \L \, Y_j. 
%\]
%%
%From the perturbed step equation we have
%%
%\begin{eqnarray*}
%\widetilde{y}_{n+1} &=&  y_n +   h\, \sum_{j=1}^s\, b_{j}\,   f(\widetilde{Y}_j)  +   h\, \sum_{j=1}^s\,  \widehat{b}_{j} \L \,  \widetilde{Y}_j \\
%&=&  y_{n+1} + h\, \sum_{j=1}^s\, b_{j}\,   f'_j \, \delta Y_j +   h\, \sum_{j=1}^s\,  \widehat{b}_{j} \L \, \delta Y_j 
%\end{eqnarray*}
%%
%and
%
\begin{eqnarray*}
\delta  y_{n+1} &=& h\, \sum_{j=1}^s\, b_{j}\,   f'_j \, \delta Y_j +   h\, \sum_{j=1}^s\,  \widehat{b}_{j} \L \, \delta Y_j.
\end{eqnarray*}

\paragraph{Non-stiff or moderately stiff problems.}
In the non-stiff or moderately stiff case where $\Vert \L \Vert = \mathcal{O}(1)$ we have
\[
\varepsilon \sim \mathcal{O}\left( h^{2k+2} \right) \quad \Rightarrow \quad
\delta Y \sim \mathcal{O}\left( h^{2k+2} \right) \quad \Rightarrow \quad \Vert \widetilde{y}_{n+1} -  y_{n+1} \Vert \sim \mathcal{O}\left( h^{2k+3} \right).
\]
For $k=0$ and $k=1$ correction iterations we have the following results.

\begin{theorem}
\label{thm:LIRK AMF order}
If a LIRK method of order higher than $2$ is applied to a nonstiff or moderately stiff case of \eqref{eqn:semi-linear} with the AMF technique according to \eqref{eqn:lirk-hong}, then the order of the method will reduce to second order.
\end{theorem}

\begin{theorem}
\label{thm:LIRK AMFR1 order}
If a LIRK method of order $3$ or $4$ is applied to a nonstiff or moderately stiff case of \eqref{eqn:semi-linear} with the AMF technique according to \eqref{eqn:lirk-hong}, 
and one correction iteration \eqref{eqn:one-ieration}
%
%\begin{eqnarray}
%\label{eqn:correction1}
%\left(\I - h\,\gamma\,\widetilde{\L} \right) Y_i^{(0)} &=&  \ell_i \\
%\label{eqn:correction2}
% Y_i^{(1)}  & =&  Y_i^{(0)} - \left( \I - h \, \gamma \, \widetilde{\L} \right)^{-1} \cdot \left( \left( \I - h \, \gamma \, \L \right)\, Y_i^{(0)} - \ell_i \right),
%\end{eqnarray}
%
is applied to each stage value $Y_i$, the order of the method is the same as that of the original method.
\end{theorem}
\begin{remark}
 Taking more iterations in the correction procedure may further reduce the linear system solution errors.
\end{remark}

\begin{remark}
 Correspondingly two iterations of the correction procedure are needed for a LIRK method of order $5$ or $6$ since $k=2$ yields an error in the solution of magnitude $\mathcal{O}\left( h^7 \right)$.  The idea can be extended to arbitrarily high order methods.
\end{remark}

\paragraph{Very stiff problems.}
In the highly stiff case a more complex analysis based on \eqref{eqn:delta-stages-2} is called for since $\Vert h\L \Vert$ can be large and $\Vert h\L \cdot \delta Y_j \Vert$ can also become large. We consider LIRK methods with a stiffly accurate implicit component $\widehat{b}_i = \widehat{a}_{s,i}$. We have that
\begin{eqnarray*}
Y_s  &=& y_n +   h\, \sum_{j=1}^{i-1}\, a_{s,j}\,   f({Y}_j)  +  h\, \sum_{j=1}^{s}\,  \widehat{a}_{s,j} \,  \L  \,{Y}_j, \\ 
y_{n+1} 
\iflong
&=&  y_n +   h\, \sum_{j=1}^s\, b_{j}\,   f(Y_j)  +   h\, \sum_{j=1}^s\,  \widehat{b}_{j} \L \, Y_j \\
\fi
&=& Y_s + h\, \sum_{j=1}^{s}\, (b_j-a_{s,j})\,   f({Y}_j). 
\end{eqnarray*}
The corresponding error equation
\begin{eqnarray*}
\delta y_{n+1}  &=& \delta Y_s + h\, \sum_{j=1}^{s}\, (b_j-a_{s,j})\,   f'_j\, \delta{Y}_j 
\end{eqnarray*}
and the non-stiff condition \eqref{eqn:nonstiff-f-jacobian} reveal that
the error in the solution is of the same size as the error in the linear solvers
\[
 \Vert \widetilde{y}_{n+1} -  y_{n+1} \Vert \sim \Vert \delta Y\Vert \sim \Vert \varepsilon \Vert.
\]

%%%%%%%%%%%%%%%%%%%%%%%%%%%%%
\subsection{Stability considerations}
%%%%%%%%%%%%%%%%%%%%%%%%%%%%%

Following Calvo and Gerisch \cite{Calvo_2005_IMEX_AMF}  we perform stability analysis using the following scalar test problem:
\begin{equation}
\label{eqn:test-stability}
h f(y) = i w y,  
\quad h\L = z = z_1 + z_2, \quad h\widetilde{\L} =  z_1 + z_2 - \gamma z_1 z_2.
\end{equation}
The LIRK method \eqref{eqn:lirk} applied to the test problem \eqref{eqn:test-stability} gives:
\begin{eqnarray*}
y_{n+1} &=& R(z,iw)\,  y_n, \\
R(z,iw) &=& 1 +  (iw\, b + z \widehat{b})^T\, \left( \I - z \widehat{A} - iw A \right)^{-1}\,\mathbf{1},\\
R(\infty,iw) &=& 1 - \widehat{b}^T\, \widehat{A}^{-1}\,\mathbf{1}.
\end{eqnarray*}
The LIRK+AMF method \eqref{eqn:lirk-hong} applied to the test problem \eqref{eqn:test-stability}  gives:
\iflong
\begin{eqnarray*}
(1-\gamma z_1)(1-\gamma z_2)\, \widetilde{Y}_i  &=& y_n +   i w\, \sum_{j=1}^{i-1}\, a_{i,j}\,   \widetilde{Y}_j  +  (z_1+z_2)\, \sum_{j=1}^{i-1}\,  \widehat{a}_{i,j} \, \widetilde{Y}_j, \\
(1+\gamma^2 z_1 z_2)\, \widetilde{Y}_i  &=& y_n +   i w\, \sum_{j=1}^{i-1}\, a_{i,j}\,   \widetilde{Y}_j  +  (z_1+z_2)\, \sum_{j=1}^{i}\,  \widehat{a}_{i,j} \, \widetilde{Y}_j, \\
\end{eqnarray*}
and therefore
\fi
\begin{eqnarray*}
\widetilde{Y}_i  &=& \left( (1+\gamma^2 z_1 z_2) \I - z \widehat{A} - iw A \right)^{-1}\, \mathbf{1}\, y_n, \\
\widetilde{y}_{n+1} &=&  y_n +   i w\, \sum_{j=1}^s\, b_{j}\,  \widetilde{Y}_j  +   (z_1+z_2)\, \sum_{j=1}^s\,  \widehat{b}_{j}  \,  \widetilde{Y}_j. 
\end{eqnarray*}
Consequently, for very stiff linear components the overall scheme is weakly stable:
\begin{eqnarray*}
y_{n+1} &=& R(z_1,z_2,iw)\,  y_n, \\
R(z_1,z_2,iw) &=& 1 +  (iw\, b + z \widehat{b})^T\, \left( (1+\gamma^2 z_1 z_2)\I - z \widehat{A} - iw A \right)^{-1}\,\mathbf{1},\\
R(\infty,\infty,iw) &=& 1.
\end{eqnarray*}

For the scheme with one refinement step we have:
\iflong
\begin{eqnarray*}
\left( 1-\gamma z_1 \right)\left( 1-\gamma z_2 \right)\, {Y}_i^{(1)}  &=& 1 +   iw\, \sum_{j=1}^{i-1}\, a_{i,j}\,   \widetilde{Y}_j  +  z\, \sum_{j=1}^{i-1}\,  \widehat{a}_{i,j} \,   \, \widetilde{Y}_j, \\
\left( 1+\gamma^2 z_1 z_2 \right)\, {Y}_i^{(1)}   &=& 1 +   iw\, \sum_{j=1}^{i-1}\, a_{i,j}\,   \widetilde{Y}_j  +  z\, \sum_{j=1}^{i}\,  \widehat{a}_{i,j} \,   \, \widetilde{Y}_j, \\
\left( 1-\gamma z_1 \right)\left( 1-\gamma z_2 \right)\, \widetilde{Y}_i &=& 
 \left( 1  - \gamma \,(z_1+z_2) + 2 \, \gamma^2 z_1 z_2 \right)\, {Y}_i^{(1)},
\end{eqnarray*}
and therefore
\fi
\begin{eqnarray*}
 \tau &=& \frac{1  - \gamma \,(z_1+z_2) + 2 \, \gamma^2 z_1 z_2}{\left( 1-\gamma z_1 \right)\left( 1-\gamma z_2 \right)\left( 1+\gamma^2 z_1 z_2 \right)}, \\
 \widetilde{Y}  &=& \left(\tau^{-1} \I - iw A - (z_1+z_2)\, \widehat{A} \right)^{-1}\;  \mathbf{1} \, y_n, \\
\widetilde{y}_{n+1} &=&  y_n + \left( iw b + z \widehat{b} \right)^T \cdot \left(\tau^{-1} \I - iw A - (z_1+z_2)\, \widehat{A} \right)^{-1}\;  \mathbf{1}\, y_n.
\end{eqnarray*}
The corresponding stability function is
\begin{eqnarray*}
y_{n+1} &=& R(z_1,z_2,iw)\,  y_n, \\
R(z_1,z_2,iw) &=& 1 + \left( iw b + z \widehat{b} \right)^T \cdot \left(\tau^{-1} \I - iw A - z\, \widehat{A} \right)^{-1}\;  \mathbf{1},\\
R(\infty,\infty,iw) &=& 1.
\end{eqnarray*}
We have that for very stiff components the scheme with one level of refinement is weakly stable. The refinement does not improve the overall stability properties of the scheme, only the accuracy. The correction step of Calvo and Gerisch \cite{Calvo_2005_IMEX_AMF} 
leads to L-stable, first order methods for very stiff problems.

%%%%%%%%%%%%%%%%%%%%%%%%%%%%%
\section{Numerical experiments} \label{sec:AMF_exp}
%%%%%%%%%%%%%%%%%%%%%%%%%%%%%

We perform numerical experiments with the following methods:
\begin{itemize}
\item \textbf{LIRK3(4):} the original LIRK methods of orders three and four, respectively, proposed by Calvo, de Frutos, and Novo \cite{Calvo_2001}. The implicit parts are L-stable and stiffly accurate; 
\item \textbf{LIRK3(4)AMF:} the LIRK methods using AMF as \eqref{eqn:lirk-hong};
\item \textbf{LIRK3(4)AMFR1:} the LIRK methods using AMF together with one iteration refinement;
\item \textbf{LIRK3(4)AMFR2:} the LIRK methods using AMF together with two iterations refinement. 
\end{itemize}
In all the experiments, the error is computed in the relative $L_2$ norm as follows
\begin{equation}
E = \frac{\| u - u_{\rm ref}\|_2}{\| u_{\rm ref}\|_2}, 
\end{equation}
where $u$ is the numerical solution at the final time, and $u_{\rm ref}$ is the reference solution at the same point obtained using Matlab's \texttt{ode15s} solver with very tight tolerances (AbsTol=RelTol=$3\times10^{-14}$).

%%%%%%%%%%%%%%%%%%%%%%%%%%
\subsection{An Allen-Cahn type problem}
%%%%%%%%%%%%%%%%%%%%%%%%%%
We use the PDE test problem of Allen-Cahn type from \cite{Calvo_2005_IMEX_AMF}:
\begin{equation}
\label{eq:allen-cahn}
u_t=\Delta u + u-u^3 + f, 
\end{equation}
where $f$ is chosen to make the exact solution of the equation be
\begin{equation*}
u(t,x,y)=e^t \sin(\pi x) \sin(\pi y). 
\end{equation*}
The spatial domain is $(x,y) \in [0,1] \times [0,1]$ and the time interval is $t \in [0,1]$ (units). 
The initial conditions and Dirichlet boundary conditions are calculated from the exact solution. 

The spatial discretization uses second order central finite differences of the Laplacian on a uniform grid of size $M\times M$
\begin{equation}
\label{eq:uniform-mesh}
 (x_i,y_i)=\left(\frac{i}{M+1},\frac{j}{M+1}\right), \quad i,j=1,\dots,M. 
\end{equation}
In our tests we consider the case $M=59$. 
The discrete solution elements $U_{ij}(t) \approx u(t,x_i,y_j)$ are ordered into a vector 
\[
z = \left( U_{11}, U_{12},\dots,U_{1M}, \dots, U_{M1}, U_{M2},\dots,U_{MM}\right)^T.
\]
The resulting ODE system can then be written into the form \eqref{eqn:semi-linear} with the discrete diffusion term being the linear part. 
The largest magnitude of the eigenvalues of the Jacobian for the diffusion term is approximately $\num{2.88e+04}$. 
The discrete Laplacian operator $\L=\L_x+\L_y$ is split into two matrices corresponding to derivatives along $x$ and $y$ directions, respectively:
\begin{equation} \label{eqn:reduction}
 \L_x = \mathbf{D}_M \otimes \I_M, ~\L_y =  \I_M \otimes \mathbf{D}_M, ~ \mathbf{D}_M = \left(
 \begin{array}{c c c c c}
  -2 & 1      &        &        & \\
  1  & -2     & 1      &        & \\
     & \ddots & \ddots & \ddots & \\
     &        & 1      & -2     & 1 \\
     &        &        & 1      & -2
 \end{array}
 \right),
\end{equation}
where the symbol $\otimes$ denotes the tensor product and $\I_M$ is an identity matrix of dimension $M\times M$. 

% Figure () shows the structure of the three matrices. 
% It can be seen that $\L_1$ and $\L_2$ have smaller bandwidth than $\L$ thus the corresponding linear systems can be better exploited for acceleration.  

\begin{figure}[!htb]
\centering{
\subfigure[Temporal errors vs. number of steps ]{
  \includegraphics[width=0.95\textwidth]{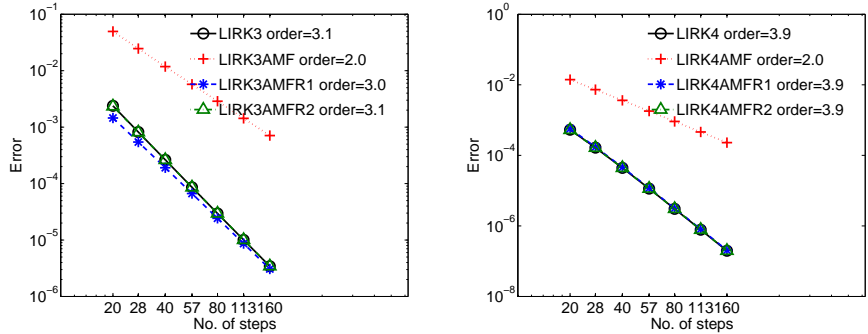}
}
\subfigure[Temporal errors vs. CPU time]{
  \includegraphics[width=0.95\textwidth]{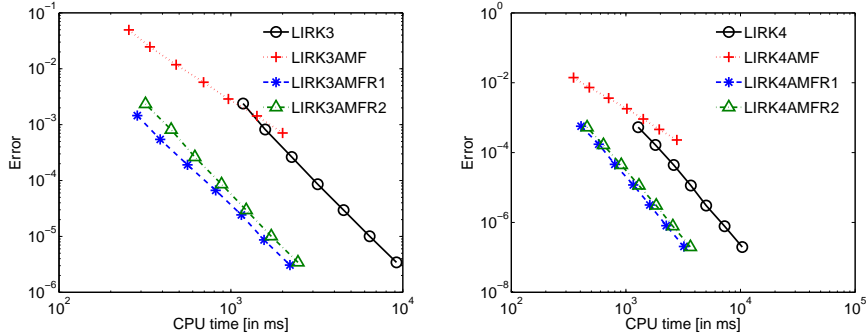}
}
}
\caption{Results for the 2D Allen-Cahn type problem \eqref{eq:allen-cahn}.}
\label{fig:allen-cahn}
\end{figure}

Figure \ref{fig:allen-cahn}(a) plots the convergence results for all the methods tested. 
As expected, both LIRK3AMF and LIRK4AMF show second order, and give less accurate results for the same time step than the underlying LIRK methods.  
All the LIRK methods with AMF and refinement perform equally well as the original LIRK methods; LIRK3AMFR1 produces slightly better results, and the full order of the underlying LIRK methods has been recovered. 

Figure \ref{fig:allen-cahn}(b) shows the corresponding work-precision diagrams. 
LIRK methods with AMF are not very competitive in terms of efficiency. 
One refinement iteration improves LIRK3AMF and LIRK4AMF significantly. 
LIRK3AMFR1 and LIRK4AMFR1 are clearly the most efficient methods among the methods of the same order. 
A second iteration does not improve accuracy, but spends compute time, and makes LIRK3(4)AMFR2 schemes slightly less efficient.
Calvo and Gerisch's approach \cite{Calvo_2005_IMEX_AMF} can achieve second order with AMF, and may only be attractive for low accuracy requirements. 
The approach with stage refinement proposed herein is competitive at all accuracy levels due to the recovery of the full order and the addition of the relatively cheap refinement procedure.
It should be noted that the separation of $\L$ into $\L_x$ and $\L_y$ in \eqref{eqn:reduction} allows for a reduction of the large linear systems with $2M$ small ones of dimension $M\times M$, 
which can lead to considerable savings in the computational cost. 
The trick is not used with this example simply for the sake of enabling a fair comparison with the results in \cite{Calvo_2005_IMEX_AMF}.

%%%%%%%%%%%%%%%%%%%%%%%%%
\subsection{Brusselator problem}
%%%%%%%%%%%%%%%%%%%%%%%%%
We next consider the two-dimensional Brusselator reaction-diffusion equation \cite[Sec. IV.10]{Hairer_book_I}
\begin{subequations}
\label{eq:brus}
\begin{eqnarray}
 u_t &=& 1 + u^2 v - (B+1)\, u + \alpha \Delta u,  \\
 v_t &=& B\, u - u^2 v + \alpha \Delta v, 
\end{eqnarray}
\end{subequations}
where $(x,y)\in[0,1]^2$, $t\in[0,1.5]$, with the Neumann boundary conditions
\[
 \frac{\partial u}{\partial \mathbf{n}} = 0, \quad \frac{\partial v}{\partial \mathbf{n}} = 0.
\]
The problem is discretized with a second order central finite difference scheme on a uniform mesh \eqref{eq:uniform-mesh}.
The stiffness of this problem increases with the value of $\alpha$ and number of grid points $M$.
We test LIRK+AMF methods with or without iterative corrections on two different cases.

\paragraph{Case 1.} 
A nonstiff stiff case described by $\alpha=0.001$, $B=3$ 
and the initial conditions
\[
 u(x,y,0)= 0.5+y, \quad v(x,y,0)= 1+5 x.
\]
with $M=39$ grid points used in each dimension. 
This gives an ODE system is of dimension $N= 3,042$.

\paragraph{Case 2.} 
A stiff case, in which the settings follow \cite[Sec. IV.10]{Hairer_book_I} where $\alpha=0.1$, $B=3.4$, and 
the initial conditions are defined as
\[
 u(x,y,0)= 22 y (1-y)^{3/2}, \quad v(x,y,0)= 22 x (1-x)^{3/2}.
\]
We choose $M=127$ for the three-way splitting test and $M=199$ for the two-way splitting test so that the resulting ODE systems are of size $N= 31,752$ and $N= 79,202$ respectively. 
The time varying Jacobian for the reaction term used in the three-way splitting test makes the linearization challenging especially when the problem is very stiff. 
For a large $M$, e.g. $M=150$, the error caused by the linearization leads to failure in solving the linear systems with direct methods after a few time steps. 
So we select a relatively smaller value of $M$ for the three-way splitting test. 
The details on the splitting setup are given later in this section. 

Table \ref{table:eigs} shows the dominant eigenvalues for these test cases which shed light on the degree of stiffness. 
Note that the Jacobian matrix for the reaction term has complex eigenvalues which makes this test problem more challenging than the previous one. 

After spatial discretization the PDE is turned into a semi-linear ODE system of the form
\[
 z' = \underbrace{\L_{dif}\,z}_{\rm diffusion} + \underbrace{\L_{rea}(t)\,z}_{\rm reaction}+R,
\]
where $z$ is the combined vector for the variables $u$ and $v$, $\L_{dif}$ and $\L_{rea}(t)$ are the Jacobian of the diffusion term and reaction term respectively, and $R$ is the rest of the terms such as boundary treatment. 
The diffusion term is stiff while the reaction term is nonstiff.  
We first apply LIRK methods with the diffusion treated implicitly and the other terms explicitly. 
Next we include the Jacobian of the reaction terms in the linear part and apply a three-way splitting strategy. 
The LU decompositions are performed per time step using sparse Gaussian elimination in MATLAB. 

\begin{table*}[h]
\small
\caption{The dominant eigenvalues (largest in magnitude) of each component.}
\centering{
\begin{tabular}{l l l l l }
\toprule
 Case    & $L_x$ & $L_y$ & $L_x+L_y$ & $L_{rea}(t_0)$ \\
\midrule
1 ($M=39$) &  \num{1.28E+1}  & \num{1.28E+1}   &  \num{2.56E+1}  & \num{1.05E+1} \\
2 ($M=127$)&  \num{6.55E+03} & \num{6.55E+03}  &  \num{1.31E+04} & \num{2.01E+1} \\
3 ($M=199$)&  \num{1.60E+04} & \num{1.60E+04}  &  \num{3.20E+04} & \num{2.01E+1} \\
\bottomrule
\end{tabular}
}
\label{table:eigs}
\end{table*} 

\paragraph{Two-way splitting.} 
A directional splitting is applied to the diffusion term which is written as the sum of derivatives in the x-direction and y-direction, $\L_{dif} = \L_x + \L_y$. 
See \eqref{eqn:reduction} for the structure of $\L_x$ and $L_y$. 
This splitting allows to reduce the linear algebra effort to solving $2M$ tridiagonal systems of dimension $M$,
all of which share the same matrix $\I- h \gamma \, \mathbf{D_M}$, and use reordered right-hand sides. 

\begin{figure}[!htb]
\centering{
\subfigure[Temporal error vs. number of steps]{
  \includegraphics[width=0.95\textwidth]{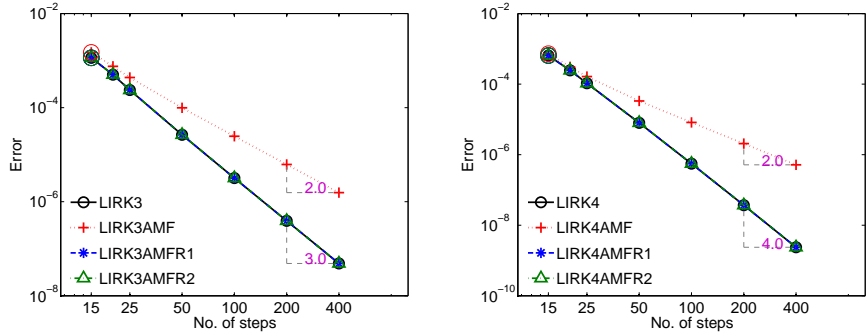}
}
\subfigure[Temporal error vs. CPU time]{
  \includegraphics[width=0.95\textwidth]{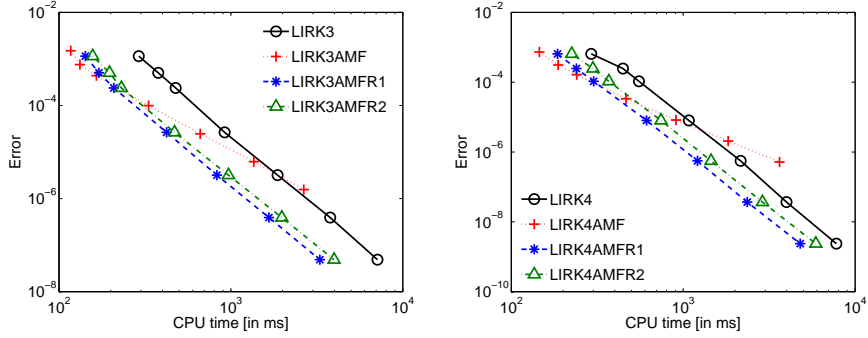}
}
}
\caption{Results for the 2D Brusselator system\eqref{eq:brus}, Case 1, with $M=39$. 
AMF is applied with a two-way splitting of the Jacobian. 
$15,20,25,50,100,200,400$ equal steps are used for the time integration of the system on the interval $[0,1]$.  
The left-most points (highlighted by a circle) on each curve indicates the maximal allowable time steps.}
\label{fig:bruss2D-case1}
\end{figure}

Figure \ref{fig:bruss2D-case1}(a) shows the convergence plots of different methods for Case 1. Both LIRK3AMF and LIRK4AMF give second order. 
With one refinement iteration the results become as accurate as those of the original LIRK methods. 
The largest allowable time steps are almost the same for all methods, implying good stability properties of LIRK methods with AMF. 
Figure \ref{fig:bruss2D-case1}(b) presents the corresponding work-precision diagrams. 
It can be seen that LIRK methods with one refinement iteration yield the best efficiency. 
To achieve the same accuracy level, LIRK3AMFR1 is about $2.2$ times faster than LIRK3 and LIRK4AMFR1 is about $1.6$ times faster than LIRK4.   

\begin{figure}[!htb]
\centering{
\subfigure[Temporal error vs. number of steps]{
  \includegraphics[width=0.95\textwidth]{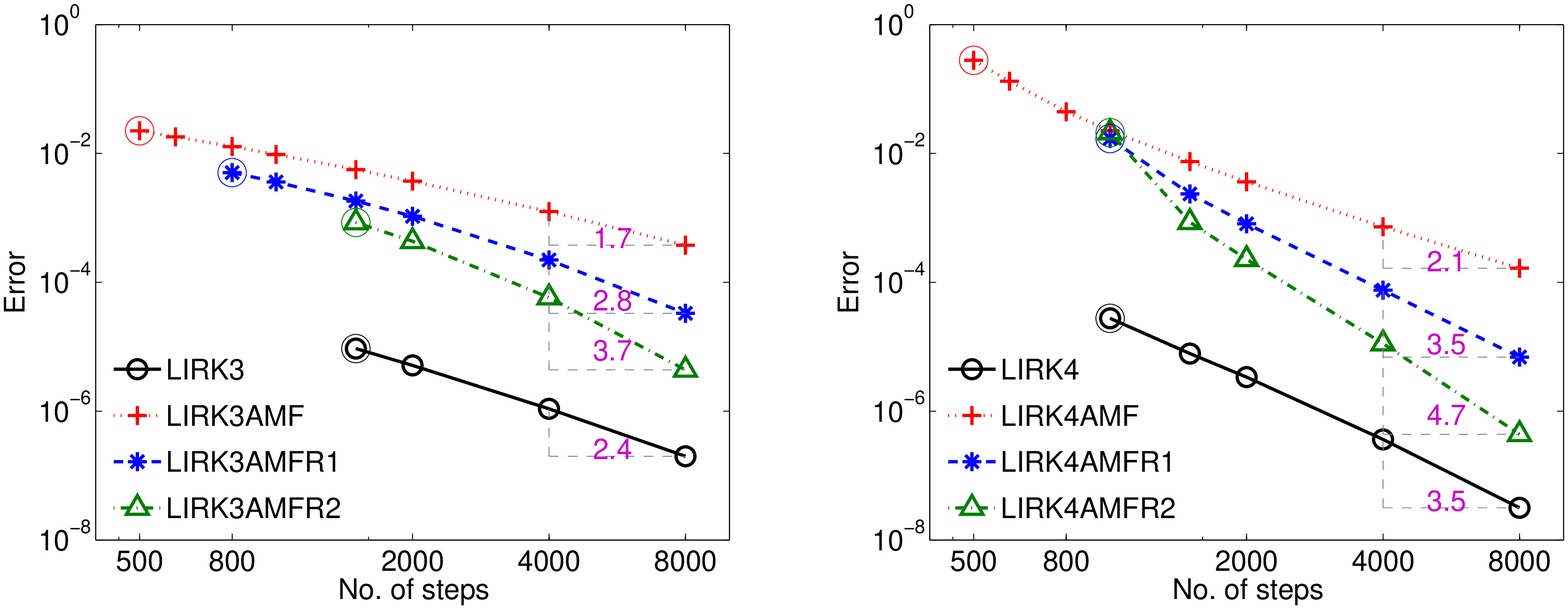}
}
\subfigure[Temporal error vs CPU time]{
  \includegraphics[width=0.95\textwidth]{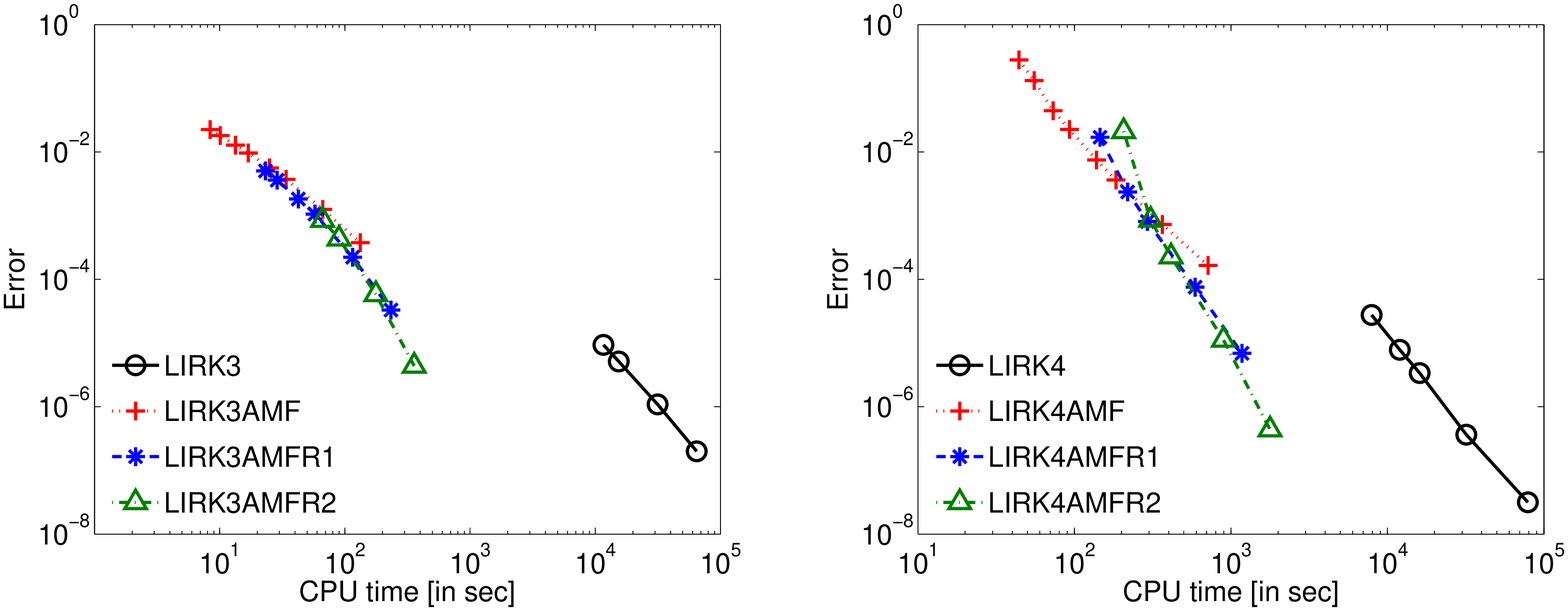}
}
}
\caption{Results for the 2D Brusselator system \eqref{eq:brus}, Case 2, with $M=199$. 
AMF is applied with a two-way splitting of the Jacobian. 
400, 500, 600, 800, 1000, 1500, 2000, 4000, 8000 equal steps are used for the time integration of the system on the interval $[0,1]$. 
The left-most points (highlighted by a circle) on each curve indicate the maximal allowable time steps.}
\label{fig:bruss2D-case2}
\end{figure}

Figure \ref{fig:bruss2D-case2} shows the convergence and work-precision diagrams for different methods for the large scale stiff case with $M=199$. 
Generally LIRK methods with AMF and two refinement iterations give more accurate results than those with AMF and one refinement iteration. 
The refinement works well and recovers the theoretical orders of the corresponding LIRK methods. 
The errors for methods with two refinement iterations approach the LIRK results at a faster rate than methods with just one refinement iteration.  
This differs from the results of the first case, but is in line with the theoretical prediction. 
Another notable advantage of the AMF technique is the gain in term of stability. 
If no refinement is employed, the maximal time step size can be at least two times larger than that allowed by the original LIRK methods for both third-order and fourth-order schemes. 
However, the gain disappears or shrinks when the refinement procedure is added.

In the efficiency comparison, LIRK methods with AMF and two refinement are the most effective for solutions more accurate than approximately $10^{-4}$. 
LIRK methods with AMF provide a good compromise between accuracy and speed since they can use a maximal allowable time step size that is at least two times larger than LIRK methods,  and run significantly faster than LIRK methods.  

\paragraph{Three-way splitting.} 
We use a three-way splitting of the linear part
\[
 \L = \L_x + \L_y + \L_{rea}(t)
\]
where the Jacobian of the reaction terms is treated implicitly.
$\L_{rea}(t)$ contains four blocks, each of which is a diagonal matrix. It is updated at each time step. 

The linear system associated with $\I- h \gamma \, \L_{rea}$ can be reduced to $M^2$ smaller systems of dimension two which can be solved separately at each grid point. We choose to solve the large sparse system directly as an explicit decoupling does not
lead to a clear performance gain for a serial implementation. 

\begin{figure}[!htb]
\centering{
\subfigure[Temporal error vs. number of steps]{
  \includegraphics[width=0.95\textwidth]{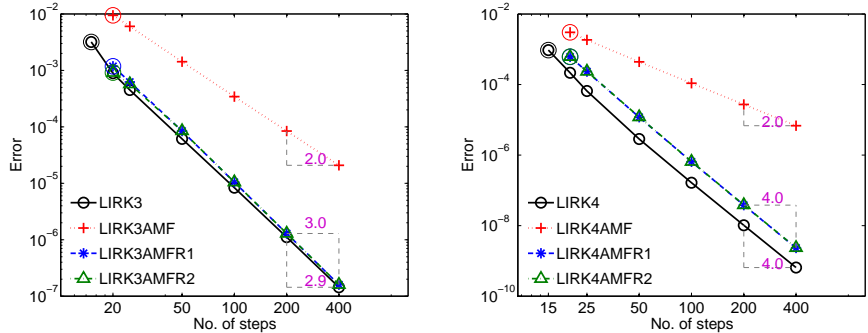}
}
\subfigure[Temporal error vs. CPU time]{
  \includegraphics[width=0.95\textwidth]{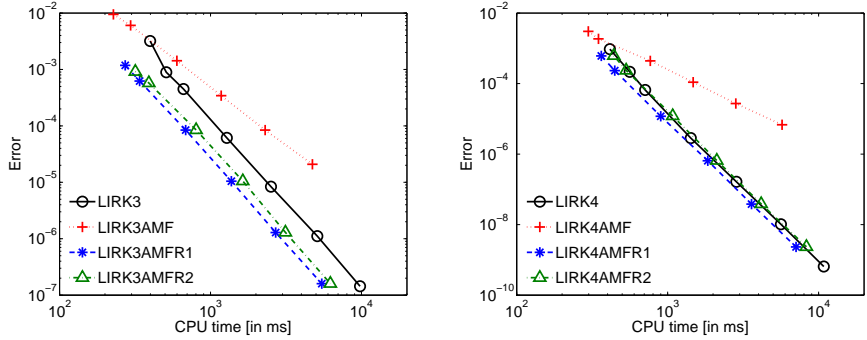}
}
}
\caption{Results for the 2D Brusselator system \eqref{eq:brus}, Case 1, with $M=39$. 
AMF is applied with a \textit{three-way} splitting of the Jacobian. 
15, 20, 25, 50, 100, 200, 400 equal steps are used for the time integration of the system on the interval $[0,1]$.  
The left-most points (highlighted by a circle) on each curve indicates the maximal allowable time steps. The numbers inside the triangle give the convergence order.}
\label{fig:bruss2D-case1-3w}
\end{figure}

\begin{figure}
\centering{
  \subfigure[Temporal error vs. number of steps]{
  \includegraphics[width=0.95\textwidth]{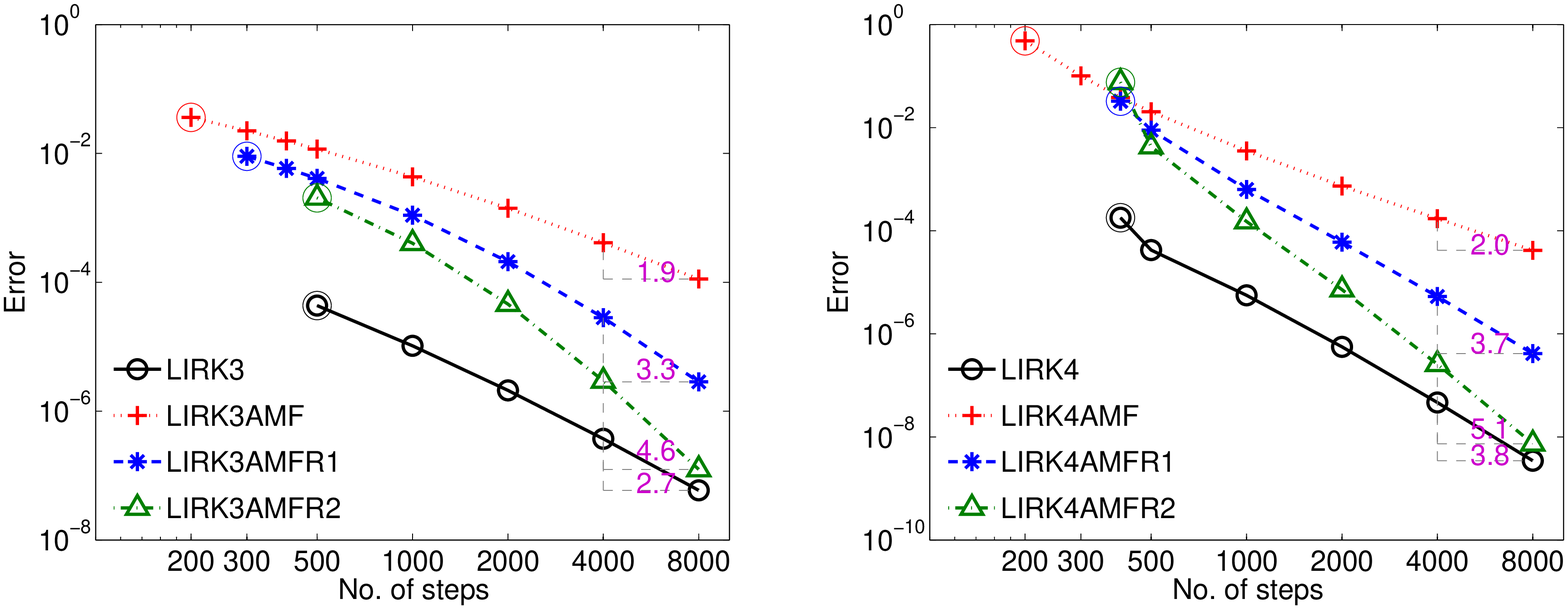}
  }
  \subfigure[Temporal error vs. CPU time]{
  \includegraphics[width=0.95\textwidth]{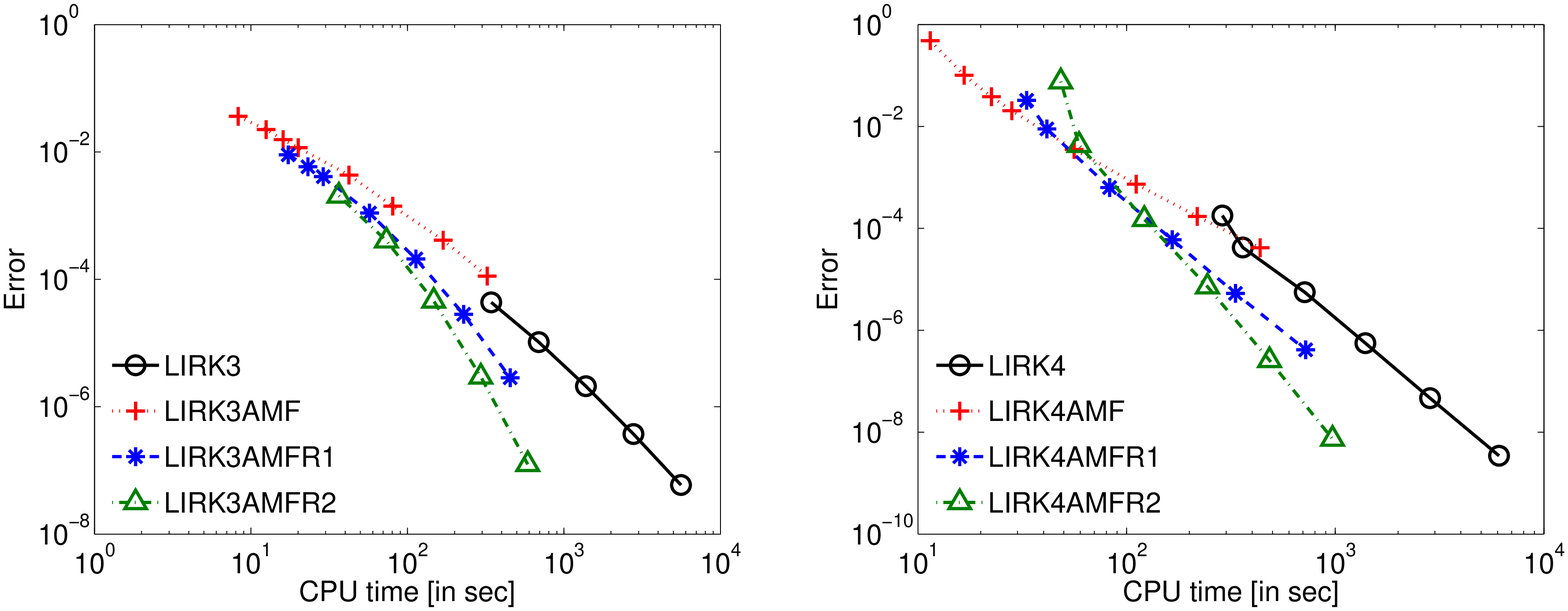}
  }
}
\caption{Results for the 2D Brusselator system \eqref{eq:brus}, Case 2, with $M=127$. 
AMF is applied with a three-way splitting of the Jacobian. 
100, 200, 300, 400, 500, 1000, 2000, 4000, 8000 equal steps are used for the time integration of the system on the interval $[0,1]$. 
The left-most points (highlighted by a circle) on each curve indicates the maximal allowable time steps. }
\label{fig:bruss2D-case2-3w}
\end{figure}

The results are shown in Figure \ref{fig:bruss2D-case1-3w}, and \ref{fig:bruss2D-case2-3w}. 
For Case 1, the refinement procedure can successfully improve the order from $2$ to the theoretical order. 
AMF without refinement is not competitive in terms of efficiency, 
but AMF with refinement yields some performance gain for third-order schemes and comparable results with LIRK methods for fourth-order schemes. 
This is due to the fact that the system is relatively small and additional cost is brought in to solve the linear system associated with $\L_{rea}(t)$. 

For the stiff case, the Jacobian for the implicit part $\L$ makes the linear system difficult to solve with direct methods. 
One LU decomposition of the system may take over $5000$ seconds. 
To improve the performance of LIRK methods, we make use of the reordering algorithm \texttt{symamd} in MATLAB before solving the linear systems. 
The reordering can also help reduce the bandwidth of the sparse matrices as similar to purpose of the splitting schemes we used.  

The convergence orders for the stiff case are very close to the two-way splitting test results. 
In terms of efficiency, there is still considerable performance gain for AMF with refinement, especially for AMF with two refinement iterations. 
The savings in CPU time may mainly come from reducing the big system into multiple small systems, which is another advantage of the application of AMF. 


\section{Conclusions} \label{sec:AMF_con}

We have applied approximate matrix factorization to high order linearly implicit Runge-Kutta methods for solving semi-linear systems of differential equations. 
The factorization (splitting) error brought by AMF leads to severe order degradation, especially for high order Runge-Kutta methods. 
The existing approach to recover second order is based on correction applied to the next step solution \cite{Calvo_2005_IMEX_AMF}. In this work  the full order of the 
underlying methods is recovered by correcting stage values via a refinement procedure based on the idea of simplified Newton iterations.  

We have performed error analysis for the linear system solutions with AMF, and investigated how this errors affect the next step solution.
In the non-stiff and mildly stiff case the full order of the underlying method can be recovered using a fixed, small number of refinement iterations.
In the very stiff case the number of iterations can be large since the convergence can deteriorate with increasing stiffness. 
A stability analysis reveals that the stage refinement procedure does not improve the overall stability of the LIRK+AMF method. When AMF is used the resulting schemes are only weakly stable for very stiff problems. Consequently, this application of AMF is attractive for mildly stiff problems, but may not work well for very stiff systems.

Numerical experiments on a variety of test problems of different sizes and different degrees of stiffness validate the theoretical findings on the accuracy and stability of high order linearly implicit Runge-Kutta methods when AMF is used. 
The results also show that the proposed approach can improve the efficiency of high order linearly implicit Runge-Kutta methods significantly and thus is attractive for solving large scale mildly stiff systems such as diffusion-reaction equations.  
Furthermore, our tests on the three-way splitting demonstrate that our methods can also efficiently deal with problems where stiffness comes from both diffusion and reaction terms. 
Though we considered LIRK schemes up to order four, the general framework developed herein can be applied to higher order LIRK methods, and could be extended to study the use of AMF with implicit-explicit multistep methods and implicit-explicit general linear methods.

\section*{Acknowledgements}

This work has been supported in part by NSF through awards NSF CMMI--1130667, 
NSF CCF--1218454, NSF CCF--0916493, NSF DMS--1419003,
 AFOSR FA9550--12--1--0293--DEF, AFOSR 12-2640-06,
and by the Computational Science Laboratory at Virginia Tech.

\bibliographystyle{plain}
\bibliography{zhang}

\appendix
%%%%%%%%%%%%%%%%%%%%%%%%%%%%%%%%%%%%%%%%%%%%%%%%%%%%%%%
\section{LIRK methods}\label{sec:coefficients} 
%%%%%%%%%%%%%%%%%%%%%%%%%%%%%%%%%%%%%%%%%%%%%%%%%%%%%%%
The coefficients of the LIRK3 method \cite{Calvo_2001}: 
\[
\renewcommand{\arraystretch}{1.5}
\tabcolsep=5pt
\begin{tabular}{r | l l l l }
$0$ & $0$ &  &  &  \\
$\gamma$ & $0$ & $\gamma$ &  &  \\
$\frac{1+\gamma}{2}$ & $0$ & $\frac{1-\gamma}{2}$ & $\gamma$ &  \\
$1$ & $0$ & $b_2$ & $b_3$ & $\gamma$\\
 \hline
 & $0$ & $b_2$ & $b_3$ & $\gamma$ \\
\end{tabular}\quad \quad \quad \quad  \quad 
\begin{tabular}{r | l l l l }
$0$ & $0$ &  &  &  \\
$\gamma$ & $\gamma$ & $0$ &  &  \\
$\frac{1+\gamma}{2}$ & $\frac{1-\gamma}{2}-a_{32}$ & $a_{32}$ & $0$ &  \\
$1$ & $0$ & $1-a_{43}$ & $a_{43}$ & $0$\\
 \hline
 & $0$ & $b_2$ & $b_3$ & $\gamma$ \\
\end{tabular},
\]
where $b_2 = -\frac{3\gamma^2}{2}+4\gamma-\frac{1}{4}$ and 
$b_3 = \frac{3\gamma^2}{2}-5\gamma+\frac{5}{4}$. 
And the choice for the free parameter is $\gamma=0.435866521508459$ and $a_{32}=0.35$.

The coefficients of the LIRK4 method \cite{Calvo_2001}: 
\[
\renewcommand{\arraystretch}{1.5}
\tabcolsep=5pt
\begin{tabular}{r | l l l l l l}
 $0$              &     &                    &                     &                  &                  & \\
 $\frac{1}{4}$    & $0$ & $\frac{1}{4}$      &                     &                  &                  & \\  
 $-\frac{3}{4}$   & $0$ & $\frac{1}{2}$      & $\frac{1}{4}$       &                  &                  & \\
 $-\frac{11}{20}$ & $0$ & $\frac{17}{50}$    & $-\frac{1}{25}$     & $\frac{1}{4}$    &                  & \\
 $-\frac{1}{2}$   & $0$ & $\frac{371}{1360}$ & $-\frac{137}{2720}$ & $\frac{15}{544}$ & $\frac{1}{4}$    & \\
 $1$              & $0$ & $\frac{25}{24}$    & $-\frac{49}{48}$    & $\frac{125}{16}$ & $-\frac{85}{12}$ & $\frac{1}{4}$ \\
 \hline
                  & $0$ & $\frac{25}{24}$    & $-\frac{49}{48}$    & $\frac{125}{16}$ & $-\frac{85}{12}$ & $\frac{1}{4}$
\end{tabular}\quad  \quad \quad \quad  \quad 
\begin{tabular}{r | l l l l l l}
 $0$              & $0$               &                 &                    &                   &                 & \\
 $\frac{1}{4}$    & $\frac{1}{4}$     & $0$             &                    &                   &                 & \\  
 $-\frac{3}{4}$   & $-\frac{1}{4}$    & $1$             & $0$                &                   &                 & \\
 $-\frac{11}{20}$ & $-\frac{13}{100}$ & $\frac{43}{75}$ & $\frac{8}{75}$     & $0$               &                 & \\
 $-\frac{1}{2}$   & $-\frac{6}{85}$   & $\frac{42}{85}$ & $\frac{179}{1360}$ & $-\frac{15}{272}$ & $0$             & \\
 $1$              & $0$               & $\frac{79}{24}$ & $-\frac{5}{8}$     & $\frac{25}{2}$    & $-\frac{85}{6}$ & $0$ \\
 \hline
                  & $0$              & $\frac{25}{24}$ & $-\frac{49}{48}$   & $\frac{125}{16}$  & $-\frac{85}{12}$ & $\frac{1}{4}$
\end{tabular}.
\]

\end{document}